\documentclass[10pt,a4paper]{article}
\usepackage{amsfonts}
\usepackage{amssymb}
\usepackage{amsmath}
\usepackage{textcomp}
\usepackage{bbm}
\usepackage{stmaryrd}

\usepackage{graphicx}
\usepackage{theorem}

\def \HHH{\mathcal H}

\def \NN{\mathbbm N}
\def \ZZ{\mathbbm Z}
\def \ZZ{\mathbbm Z}

\def \RR{\mathbbm R}

\def \DD{\mathcal D}

\def \EE{\mathbbm E}

{\theorembodyfont{\rmfamily} }
\newtheorem{theorema}{Theorem}[section]
\newtheorem{prop}{Proposition}[section]
\newtheorem{lemma}{Lemma}[section]

\newcommand{\ra}{\rightarrow}

\newcommand{\lra}{\longrightarrow}

\newcommand{\D}{\Delta}
\newcommand{\g}{\gamma}

\newcommand{\subs}{\subseteq}

\newcommand{\wt}{\widetilde}
\newcommand{\ee}{\varepsilon}

\newcommand{\phih}{\varphi}
\newcommand{\mx}{\vee}
\newcommand{\mn}{\wedge}
\newcommand{\x}{\times}

\newcommand{\s}{\sigma}
\newcommand{\ls}{\langle}
\newcommand{\rs}{\rangle}

\newcommand{\sm}{\setminus}

\newcommand{\bbar}{\overline}
\newcommand{\fr}{\frac{1}}
\newcommand{\PP}{\mathbbm{P}}

\newcommand{\T}{\mathbbm{T}}
\newcommand{\MM}{\mathbbm{M}}

\newcommand{\1}{\mathbbm{1}}

\newcommand{\y}{\upsilon}
\newcommand{\pd}{\partial}

\makeatletter
\newcommand{\argmin}{\mathop{\operator@font argmin}}
\makeatother
\makeatletter
\newcommand{\astop}{\mathop{\operator@font \ast}}
\makeatother
\begin{document}
\title{Hydrodynamic Limit of mean zero condensing Zero Range Processes with sub-critical initial profiles}
\date{}
\author{ Marios-Georgios Stamatakis\footnote{Department of Applied Mathematics, University of Crete, 714 09, Heraklion Crete, Greece.}}
\maketitle
\begin{abstract}\indent We prove the hydrodynamic limit of mean zero condensing Zero Range Processes with bounded local jump rate, for sub-critical initial profiles. The proof is based on H.T. Yau's relative entropy method and is made possible by a generalisation of the One Block Estimate.\\ \end{abstract}

\section{Introduction} 

\indent In this article we study the hydrodynamic limit of mean zero Zero Range Processes (ZRPs). ZRPs are interacting particle systems such that each particle $X$ jumps at an exponential rate $g(k)$ that depends only on the number $k$ of particles that occupy the same site as particle $X$, through some function $g:\ZZ_+\lra\RR_+$ called the local jump rate. Particles that jump change position according to a transition probability $p$. In the mean zero case studied here, $p$ is assumed to have mean zero.\\ 
\indent Since their introduction by Spitzer in 1970, ZRPs have attracted a lot of attention, one reason being that for particular choices of local jump rate functions $g$ they exhibit phase transition phenomena, via the emergence of mass condensation at densities above a critical density $\rho_c$. So since the equilibrium states of ZRPs are explicitly known as product measures, ZRPs can serve as simple prototype models for the study of condensation phenomena. See section 3 in \cite{Evans1} for a review of such applications of ZRPs. In this sense this work is a study of hydrodynamic limits in the presence of condensation via the example of ZRPs. \\
\indent Our approach is based on H.T. Yau's relative entropy method, which was originally applied to prove the hydrodynamic limit of Ginzburg-Landau models in \cite{Yau}. It has been applied to mean zero ZRPs, only under assumptions on the jump rate $g$ that force the critical density $\rho_c$ to be infinite. Under these assumptions there exists an equilibrium state $\nu_\rho$ corresponding to the density $\rho$ for each $\rho\geq 0$, and the mean local jump rate $\Phi(\rho)$ under the one-site marginal of $\nu_\rho$ is $C^\infty$ with strictly positive derivative. In this case, the hydrodynamic equation of the  ZRP is the non-linear diffusion equation 
\begin{eqnarray}\label{HL}\qquad
\pd_t\rho=\D_\Sigma\Phi(\rho)
\end{eqnarray}
where $\D_\Sigma$ is the differential operator 
$\D_\Sigma=\sum_{i,j=1}^d\s_{ij}\pd^2_{ij}$
and $\Sigma=(\s_{ij})_{1\leq i,j\leq d}$ is the covariance matrix of the step distribution $p$. This application of the relative entropy method is contained in the finite volume case in chapter 6 of \cite{Landim} for bounded initial profiles and extended to the infinite volume case in \cite{LandMour}. It relies on the existence of sufficiently smooth classical solutions of the hydrodynamic equation (\ref{HL}).\\
\indent Here we prove the hydrodynamic limit of mean zero ZRPs under weaker assumptions that do not exclude ZRPs with finite critical density $\rho_c$, when starting from a local equilibrium of some sub-critical profile $\rho_0$, i.e$.$ such that it takes values in some closed sub-interval of $(0,\rho_c)$. 
Our assumption on the sub-criticality of the initial profile can be considered as the adaptation of the boundedness assumption of the initial profile to the case of finite critical density, and it serves exactly the same purpose. It allows to obtain a sufficiently smooth solution to the hydrodynamic equation (\ref{HL}) that satisfies the maximum principle and to apply the relative entropy method.\\
\indent There remains one problem. The one-block estimate, a main tool in the application of the relative entropy method, has been proven so far for ZRPs with infinite critical density. Its generalisation to ZRPs with finite critical density given here is of interest in its own right, since the one-block estimate is a main tool in any approach to the hydrodynamic limit of ZRPs. We then have all the main ingredients required for the application of the relative entropy method at our disposal and a careful application yields the hydrodynamic limit of mean zero condensing ZRPs.\\
\indent Since classical solutions of the hydrodynamic equation starting from a sub-critical profile satisfy the maximum principle, a main consequence of this more general version of the hydrodynamic limit is that the sub-criticality of the initial profile ensures that no phase transition will occur in the hydrodynamic limit, establishing thus a maximum principle at the hydrodynamic limit, even in ZRPs with finite critical density.

\section{Preliminaries}
We give in this section the definition of ZRPs and describe their equilibrium states and the known hydrodynamic limit and its assumptions. A standard reference for the material in this section is the textbook \cite{Landim} and the article \cite{Stefan}.\\
\indent A {\it{local jump rate}} is a function $g:\ZZ_+\lra\RR_+$ such that
\begin{itemize}
\item[(a)] $g(k)=0\quad\Longleftrightarrow\quad k=0$
\item[(b)] $\|g'\|_\infty:=\sup_{k\in\ZZ_+}|g(k+1)-g(k)|<+\infty$, and
\item[(c)] $\varphi_c:=\liminf_{k\rightarrow\infty}\sqrt[k]{g!(k)}>0$,
\end{itemize}
where $g!(0):=1$ and $g!(k):=g(1)\cdot g(2)\cdot\dots\cdot g(k)$. An {\it{elementary step distribution}} is a probability distribution $p\in\PP\ZZ^d$ (where for any polish space $M$ we denote by $\PP M$ the space of all Borel probability measures on $M$) such that its support ${\rm{spt}}(p):=\{z\in\ZZ^d|p(z)>0\}$ is bounded and generates $\ZZ^d$, that is for any $z\in\ZZ^d$ there exist $m\in\NN$ and $z_1,\cdots,z_m\in{\rm{spt}}(p)$ such that $z=z_1+\cdots+z_m$. The assumption that ${\rm{spt}}(p)$ generates $\ZZ^d$ implies that the covariance matrix $\Sigma=(\s_{ij})_{1\leq i,j\leq d}=\big(\int_{\ZZ^d}k_ik_jdp(k)\big)\in\RR^{d\x d}$ 
of $p\in\PP\ZZ^d$ is strictly positive definite, that is $\ls\y,\Sigma\y\rs>0$ for all $\y\in\RR^d\sm\{0\}$.\\
\indent It is convenient to consider ZRPs that evolve on the discrete $d$-dimensional tori $\T_N^d\cong(\,^\ZZ/_{N\ZZ})^d\cong\{0,1,\dots,N-1\}^d$, $N\in\NN$, and consider the limit as $N\ra\infty$. The state space of a ZRP evolving on $\T_N^d$ is the space of configurations $\MM_N^d$ consisting of all functions $\eta:\T_N^d\lra\ZZ_+$ so that given $\eta\in\MM_N^d$, $\eta_x$ is the number of particles occupying the site $x\in\T_N^d$. We will denote by $\eta(x):\MM_N^d\lra\ZZ_+$, $x\in\T_N^d$, the natural projections. The {\it{ZRP of local jump rate function}} $g$ {\it{and elementary step distribution}} $p$ {\it{on the discrete torus}} $\T_N^d$ is the unique Markov jump process on the Skorohod space $D(\RR_+;\MM_N^d)$ with generator $L^N$ given by the formula 
$$L^Nf(\eta)=\sum_{x,y\in\T_N^d}\big\{f(\eta^{x,y})-f(\eta)\big\}g\big(\eta_x\big)p(y-x),$$
where  $\eta^{x,y}$ is the configuration resulting from $\eta$ by moving a particle from $x$ to $y$. We denote
by $S^N_t:\MM_N^d\lra\PP\MM_N^d$ the transition semi-group of the ZRP.\\
\indent The total number of particles is conserved by the dynamics of the ZRP, and by the assumption that the support of the elementary step distribution $p$ generates $\ZZ^d$, the communication classes of the ZRP are exactly the hyperplanes 
$$\MM_{N,K}^d:=\bigg\{\eta\in\MM_N^d\Big|\sum_{x\in\T_N^d}\eta(x)=K\bigg\},\quad K\in\ZZ_+,$$
with a fixed number of particles. So for each $(N,K)\in\NN\x\ZZ_+$, there exists a unique equilibrium state $\nu_{N,K}^d\in\PP\MM_N^d$ concentrated on $\MM_{N,K}^d$. We will refer to the family $\{\nu_{N,K}\}_{(N,K)\in\NN\x\ZZ_+}$ as {\it{the canonical ensemble}} of the ZRP.\\
\indent The function $Z\equiv Z_g:\RR_+\lra[1,\infty]$ defined by the power series $$Z(\varphi)\equiv Z_g(\varphi):=\sum_{k=0}^\infty\frac{\varphi^k}{g!(k)}
$$ is called the {\it{normalising partition function associated to}} the local jump rate function $g$. The radius of convergence of $Z$ is 
$\varphi_c=\liminf_{k\ra\infty}\sqrt[k]{g!(k)}$ and so assumption (c) of local jump rate functions guaranties that $Z$ has non-trivial domain of convergence. Obviously any partition function $Z:\RR_+\lra[1,+\infty]$ is $C^\infty$ on $[0,\phih_c)$ with all of its derivatives strictly positive there. By Abel's theorem on power-series, $Z$ and all of its derivatives are lower semi-continuous on $\RR_+$. For any $\phih\in\DD_Z:=\{Z<+\infty\}$, the product distribution $\bar{\nu}_{\varphi}^N\equiv\bar{\nu}_{\varphi,g}^N\in\PP\MM_N^d$ with common marginal $\bar{\nu}^1_{\varphi}\in\PP\ZZ_+$ given by 
$$
\bar{\nu}^1_\varphi\{k\}=\fr{Z(\varphi)}\frac{\varphi^k}{g!(k)},\qquad k\in\ZZ_+
$$
is called the {\it{zero range product distribution on $\T_N^d$ of rate $g$ and fugacity}} $\varphi$.\\
\indent Note that whenever $\phih_c\in\DD_Z$ the one-site zero range distribution $\bar{\nu}_{\phih_c}^1$ corresponding to the critical fugacity $\phih_c$ is defined. The zero range product distributions $\bar{\nu}_{\phih}^N\in\PP\MM_N^d$, $\phih\in \DD_Z$, are equilibrium distributions, i.e$.$ $\bar{\nu}_{\phih}^NL^N = 0$,
and translation invariant, that is 
$\tau_{x*}\bar{\nu}_{\phih}^N=\bar{\nu}_{\phih}^N$ for all $x\in\T_N^d$, where $\tau_x:\MM_N^d\lra\MM_N^d$ denotes the translation operator $\tau_x\eta(y)=\eta(x+y)$. In fact they are the only translation invariant equilibrium states of the ZRP that are also product measures.\\
\indent As is well known, the zero range product distributions can be reparametrised by  the density. The {\it{mean density function}} $R:\DD_Z\lra[0,+\infty]$ defined by 
\begin{eqnarray}\label{MeanDensityMainFormula}\qquad
R(\phih)=E_{\bar{\nu}_{\phih}^N}[\eta(0)]=\int kd\bar{\nu}^1_\phih(k)=\frac{\phih Z'(\phih)}{Z(\phih)}
\end{eqnarray} is continuous on $\DD_Z$ and $C^\infty$ and strictly increasing on $[0,\phih_c)$ (\cite{Landim}). Consequently, its inverse $\Phi:=R^{-1}$ is well defined on $R(\DD_Z)\supseteq[0,\rho_c)$, where 
$$\rho_c\equiv R(\phih_c):=\lim_{\phih\uparrow\phih_c}R(\phih)\in(0,\infty],$$
and $\rho_c\in R(\DD_Z)$ if and only if (iff) $\phih_c\in\DD_Z$. By reparametrising the zero-range distributions by the mean jump rate $\Phi$ we get for any $\rho\in R(\DD_Z)$ an equilibrium distribution $\nu^N_\rho$ of mean density $\rho$:
\begin{eqnarray}\label{GrandCanEns}\qquad
\nu^N_\rho:=\bar{\nu}^N_{\Phi(\rho)},\quad\rho\in R(\DD_Z).
\end{eqnarray}
We will refer to the family defined in (\ref{GrandCanEns}) as the {\it{grand canonical ensemble}} of the ZRP.\\
We note that the mean jump rate function $\Phi$ is Lipschitz with Lipschitz norm $\leq\|g'\|_\infty$ and is {\it{the mean jump rate function}} since for all $\rho\in[0,\rho_c)$ we have that
$$\Phi(\rho)=E_{\nu_\rho^N}[g(\eta(0))]=\int g(k)d\nu_\rho^1(k).$$
\indent The various cases of the set $R(\DD_Z)$ are as follows. As is proved in \cite{Landim}, whenever $\phih_c\notin\DD_Z$, that is whenever $\phih_c=+\infty$ or $\phih_c<+\infty$ and $Z(\phih_c)=+\infty$, we have that $\rho_c=+\infty$. So in this case $R(\DD_Z)\supseteq\RR_+$ and the mean jump rate function $\Phi$ is defined on all of $\RR_+$. On the other hand if $\phih_c\in\DD_Z$ then $R(\DD_Z)=[0,\rho_c]$ and in this case, as is shown by (\ref{MeanDensityMainFormula}), the critical density is infinite if $Z'(\phih_c)\equiv\sup_{\phih<\phih_c}Z '(\phih)=+\infty$ and finite if $Z'(\phih_c)<+\infty$. In particular, whenever $\rho_c<+\infty$ we have that $\phih_c\in\DD_Z$ and so the grand canonical ensemble contains an equilibrium distribution corresponding to the critical density $\rho_c$.
 To our knowledge, so far the one block estimate and the hydrodynamic limit of ZRPs have only been only under the assumption $\phih_c\notin\DD_Z$, which excludes ZRPs with finite critical density, and the aim of this article is to remove this assumption.\\

\noindent{\bf{An example with finite critical density:}} In \cite{Evans2} Evans introduces ZRPs with local jump rate function 
\begin{eqnarray}\label{EvansJumpRate}\qquad
g_b(k)=\1_{\{k\geq 1\}}\Big(1+\frac{b}{k}\Big),\quad b\geq 0.
\end{eqnarray}
It is well known (\cite{Stefan}) that $\phih_c=1$ for all $b\geq 0$, $\phih_c\notin\DD_Z$ iff $b\in[0,1]$ and that for $b>2$, the first moment of the grand canonical distribution $\nu_{\phih_c}^1$ is finite, thus leading to a finite critical density $\rho_c<\infty$.\\

\indent The explicit formula of the grand canonical ensemble permits an easy computation of the canonical ensemble, since for all $(N,K)\in\NN\x\ZZ_+$, $\rho\in[0,\rho_c)$ we have $\nu_{N,K}(\cdot)=\nu_\rho^N\{\,\cdot\,|\MM_{N,K}^d\}$. It follows that the canonical equilibrium distributions are translation invariant and so under $\nu_{N,K}$ each site has a mean number of particles equal to $K/N^d$.\newpage
\indent The phase transition in ZRPs with finite critical density has been described in \cite{Evans2} and proved rigorously in \cite{Stefan} as a continuous phase transition in the thermodynamic limit:
\begin{theorema}{\rm{(Equivalence of Ensembles)}}\label{EoE}
Let $\{\nu_{N,K}\}_{K\in\ZZ_+}$ and $\{\nu_\rho^N\}_{\rho\in R(\DD_Z)}$ be the canonical and grand canonical ensemble of the ZRP, let $\pi^L:\MM_N^d\lra\MM_L^d$, $N\geq L$, be the natural projections and set $\nu_{N,K}^L:=\pi^L_*\nu_{N,K}$. Then for fixed $L\in\NN$, for all $\rho\geq 0$ it holds that $$\lim_{\substack{N,K\ra+\infty\\K/N^d\ra\rho}}\HHH(\nu_{N,K}^L|\nu_{\rho\mn\rho_c}^L)=0.$$ 
In particular $\nu_{N,K}^L\lra\nu_{\rho\mn\rho_c}^L$ weakly as $N,K\ra\infty$ and $K/N^d\ra\rho$.
\end{theorema}

In appendix 1, corollary 1.7 in \cite{Landim} a different version of the equivalence of ensembles is proved under the additional assumption that $Z(\phih_c)=+\infty$: For each $\rho_0<+\infty$, for all cylinder functions (i.e$.$ functions that depend on a finite number of coordinates) with finite second moment with respect to the measures $\nu_{\rho}^\infty$, $\rho\in[0,\rho_0]$, it holds that $$\int fd\nu_{N,K}\lra\int fd\nu_{\rho}^\infty\qquad\mbox{as }\;N,K\ra\infty\;\;\mbox{and }\;K/N^d\ra\rho$$
uniformly over all $\rho\in[0,\rho_0]$, where $\nu_\rho^\infty:=(\nu_{\rho}^1)^{\otimes\ZZ^d}\in\PP\MM_\infty^d:=\PP\ZZ_+^{\ZZ^d}$. Of course this cannot be true for $\rho>\rho_c$ if $\rho_c<+\infty$ since even for the linear cylinder function $\eta(0)$ 
$$\int\eta(0)d\nu_{N,K}\lra\rho>\rho_c\qquad\mbox{as }\;N,K\ra\infty\;\;\mbox{and }\;K/N^d\ra\rho.$$
In other words, at the thermodynamic limit we have a mean total loss of mass equal to $\rho-\rho_c$ at each site. As it has been proven, in many cases the excess mass of all the sites is concentrated on a single random site. We refer to \cite{Stefan,Mixalhs,MixInesSuperCrit} for a detailed description of the phase separation in the Evans model.\\
\indent If the local jump rate $g$ is bounded, then by the equivalence of ensembles we get that 
$$\lim_{\substack{N,K\ra\infty\\K/N^d\ra\infty}}\int g(\eta(0))d\nu_{N,K}=\int g(\eta(0))d\nu_{\rho\mn\rho_c}^\infty=\Phi(\rho\mn\rho_c),$$ 
for all $\rho\geq 0$. As noted in \cite{Stefan}, this shows that for bounded local jump rate functions $g$ the mean jump rate function $\Phi$ should be extended on all of $\RR_+$ by \begin{eqnarray}\label{MJRExtension}
\Phi(\rho)\equiv\Phi(\rho\mn\rho_c),\qquad\mbox{for all }\rho\geq 0.
\end{eqnarray}
It turns out that this choice of $\Phi$ is the right one in order to extend the one-block estimate to ZRPs with finite critical density.\\
\indent We review next the concepts of local equilibrium and hydrodynamic limits in the context of ZRPs. For any bounded cylinder function $\Psi:\MM_\infty^d\lra\RR$ we denote by $\wt{\Psi}:\RR_+\lra\RR$ the function given by $\wt{\Psi}(\rho)=\int\Psi d\nu^\infty_{\rho\mn\rho_c}$. Let $\rho:\T^d\lra\RR_+$ be a measurable function. A sequence $\mu^N\in\PP\MM_N^d$, $N\in\NN$, is called a {\it{weak local equilibrium of profile}} $\rho$ if for all bounded cylinder functions $\Psi:\MM_\infty^d\lra\RR$, all $G\in C(\T^d)$ and all $\ee>0$ 
\begin{eqnarray}\label{WeakLocalEquil}\qquad
\lim_{N\ra\infty}\mu^N\bigg\{\bigg|\fr{N^d}\sum_{x\in\T_N^d}G\Big(\frac{x}{N}\Big)\tau_x\Psi-
\int_{\T^d}G(u)\wt{\Psi}\big(\rho(u)\big)du\bigg|>\ee\bigg\}=0.
\end{eqnarray} 
\indent Given any continuous profile $\rho:\T^d\lra\RR_+$ a particular weak local equilibrium of profile $\rho$ is the sequence $\nu^N_{\rho(\cdot)}\in\PP\MM_N^d$, $N\in\NN$, of the so-called {\it{product measures with slowly varying parameter associated to the profile}} $\rho$, i.e$.$ the sequence of the product measures with one-site marginals $\nu^1_{\rho(\frac{x}{N})}$, $x\in\T_N^d$.\\
\indent In order to prove the hydrodynamic limit of mean zero ZRPs one has to prove the {\it{conservation of weak local equilibrium in the diffusive time-scale along some measurable function}} $\rho:\RR_+\x\T^d\lra\RR_+$, i.e$.$ that starting from a weak local equilibrium $\mu_0^N\in\PP\MM_N^d$, $N\in\NN$, of some sufficiently regular initial profile $\rho_0:\T^d\lra\RR$ at time $t=0$, at each later time $t>0$ the sequence of laws $\mu_t^N:=\mu_0^NS_{tN^2}$, $N\in\NN$, of the diffusively rescaled ZRP at time $t$ is a weak local equilibrium of profile $\rho_t=\rho(t,\cdot)$, where $\rho=(\rho_t)_{t\geq 0}$ is a solution of some evolutionary type PDE, the so called {\it{hydrodynamic equation}}.\\
\indent Given a classical solution $\rho:\RR_+\x\T^d\lra\RR$ of equation (\ref{HL}), the relative entropy method controls the relative entropy of the evolution $\mu_t^N$ of the diffusively rescaled ZRP with respect to the measures with slowly varying parameter associated to the profile $\rho_t=\rho(t,\cdot)$. In effect, the relative entropy method works by proving the conservation of a slightly stronger notion of local equilibrium: A sequence $\mu^N\in\PP\MM_N^d$ is an {\it{entropy local equilibrium of profile $\rho:\T^d\lra\RR_+$}} if 
$$\lim_{N\ra\infty}\fr{N^d}H(\mu^N|\nu_{\rho(\cdot)}^N)=0.$$
\noindent According to corollary 1.3 in \cite{Landim} the notion of entropy local equilibrium is stronger than the notion of weak local equilibrium. So, the relative entropy method which proves the conservation of local equilibrium in the entropy sense, also proves the conservation of weak local equilibrium, under the slightly stronger assumption that the initial local equilibrium is an entropy local equilibrium.\\
\indent The application of the relative entropy method requires the existence of sufficiently regular classical solutions to the hydrodynamic equation. By the classic results on quasi-linear uniformly parabolic equations obtained in \cite{Lady}, the main property of quasi-linear parabolic equations required to yield existence of classical solutions to (\ref{HL}) is uniform parabolicity.  In the so far proved cases of the hydrodynamic limit of ZRPs (\cite{Landim}, \cite{LandMour}) it holds that $\rho_c=+\infty$, $\Phi$ is Lipschitz and $C^\infty$ with strictly positive derivative on $[0,\infty)$. So since the matrix $\Sigma$ is positive definite, the uniform parabolicity of (\ref{HL}) is equivalent to 
\begin{eqnarray}\label{UP}\qquad
c:=\inf_{\rho\in[0,\infty)}\Phi '(\rho)>0.
\end{eqnarray} Even under the assumption $\phih_c\notin\DD_Z$, inequality (\ref{UP}) is not always true. However, since $\Phi '$ is strictly positive, (\ref{UP}) can only fail as $\rho\ra\infty$ and so for bounded initial profiles the maximum principle permits to circumvent the loss of uniform parabolicity at large densities and still obtain sufficiently smooth classical solutions. We will describe this in more detail in the proofs section in a slightly more general case, covering the mean jump rate functions of ZRPs with finite critical density $\rho_c<+\infty$. In particular we consider the {\it{class}} $\mathcal{C}_{\rho_c}$, $\rho_c\in(0,+\infty]$, of continuous functions $\Phi:\RR_+\lra\RR_+$ with $\Phi(0)=0$ such that 
 $\Phi$ is $C^\infty$ on the interval $[0,\rho_c)$ with $\Phi '(\rho)>0$ for all $0\leq\rho<\rho_c$ and $\Phi(\rho)=\Phi(\rho\mn\rho_c)$ for all $\rho\geq 0$. A continuous function $\rho:\RR_+\x\T^d\lra\RR_+$ is a {\it{classical solution of the initial value problem for the non-linear diffusion equation with non-linearity}} $\Phi\in\mathcal{C}_{\rho_c}$ if $\rho$ is such that all differentiations in (\ref{HL}) make sense and it satisfies (\ref{HL}) with  $\rho(0,\cdot)\equiv\rho_0$.

\section{Statement of the Results}\label{Statements}

Our first result is a generalisation of the one-block estimate to condensing ZRPs.
The one-block estimate has been proven under the assumption that either (a) $\phih_c=+\infty$ or (b) $\phih_c\notin\DD_Z$ and the local jump rate $g$ has sublinear growth, which excludes ZRPs with finite critical density. Here we give the proof under the assumption that the jump rate $g$ is bounded. The proof under this assumption is useful for condensing ZRPs since it does not impose any restrictions on the behaviour of the partition function $Z$ at the critical fugacity $\phih_c<+\infty$. For each configuration $\eta\in\MM_N^d$ we denote by $\eta^\ell$ the spatial mean of $\eta$ over microscopic boxes of radius $\ell$ given by $$\eta^\ell(x)=\fr{(2\ell+1)^d}\sum_{y\in\T_N^d:|y-x|\leq\ell}\eta(y).$$
\vspace{-12pt}
 \begin{theorema}{\rm{(One-Block Estimate)}} Suppose the local jump rate $g$ of the ZRP is bounded and let $\mu_0^N\in\PP\MM_N^d$ be any sequence satisfying the $O(N^d)$-entropy assumption, i.e$.$ such 
$$C(\rho_*):=\sup_{N\in\NN}\fr{N^d}H(\mu_0^N|\nu_{\rho_*}^N)<+\infty,$$
for some (and thus for any) $\rho_*\in(0,\rho_c)$. Then
\begin{eqnarray}\label{TimeDependentOneBlockEstimate}\qquad
\lim_{\ell\ra\infty}\limsup_{N\ra\infty}\EE^N\bigg|\int_0^T\fr{N^d}\sum_{x\in\T_N^d} H\Big(t,\frac{x}{N}\Big)\Big[g(\eta_t(x))-\Phi\big(\eta^\ell_t(x)\mn\rho_c\big)\Big]dt\bigg|=0
\end{eqnarray}
for all $H\in C([0,T]\x \T^d)$, $T>0$, where $\EE^N$ denotes the expectation with respect to the diffusively accelerated law of the ZRP starting from $\mu_0^N\in\PP\MM_N^d$.
\end{theorema}

\begin{theorema}\label{MainTheorem}{\rm{(Hydrodynamic Limit)}} Suppose that the local jump rate function $g$ of the ZRP is bounded and let $\Phi$ be the mean jump rate function associated to $g$. Then any initial entropy local equilibrium $\mu_0^N\in\PP\MM_N^d$, $N\in\NN$, of profile $\rho_0\in C^{2+\theta}\big(\T^d;(0,\rho_c)\big)$ is conserved in the diffusive timescale along the unique solution $\rho:\RR_+\x\T^d\lra(0,\rho_c)$ of the initial value problem \begin{eqnarray}\label{NonLinDiff}\qquad
\begin{cases}\pd_t\rho=\D_\Sigma\Phi(\rho)\quad\mbox{in}\quad(0,\infty)\x\T^d\\\rho(0,\cdot)=\rho_0.\end{cases}
\end{eqnarray}
\indent In other words, if $H(\mu_0^N|\nu_{\rho_0(\cdot)}^N)=o(N^d)$ then $H(\mu_t^N|\nu_{\rho_t(\cdot)}^N)=o(N^d)$ for all $t>0$, where $\mu_t^N:=\mu_0^NS_{tN^2}$, and in particular (\ref{WeakLocalEquil}) holds
for all $G\in C(\T^d)$ and all $t,\ee>0$. 
\end{theorema}

\noindent\textbf{Remark:} As will be seen in the proof, one can assume the initial profile $\rho_0$ to be only of class $C(\T^d;[0,\rho_c))$, provided that the unique classical solution $\rho$ of the hydrodynamic equation (\ref{NonLinDiff}) with initial condition $\rho(0,\cdot)=\rho_0$ is such that the functions
\begin{eqnarray}\label{SufficesforGeneralInitProf}\qquad
(a)\;\;t\mapsto\log\Phi\big(m_t\big)\quad\mbox{and}\quad
(b)\;\;t\mapsto\frac{\|\D_\Sigma\Phi(\rho_t)\|_\infty+\|D^2[\Phi(\rho_t)]\|_\theta}{\Phi(m_t)} 
\end{eqnarray} belong in $L^2_{\rm{loc}}(\RR_+)$, where $m_t:=\min_{u\in\T^d}\rho_t(u)$, $\|D^2 f\|_\theta:=\max_{|\alpha|=2}|\pd^\alpha f|_\theta$, $\theta\in(0,1]$, and $|f|_\theta$ is the $\theta$-Holder semi-norm of the function $f$. By the properties of the solutions of the hydrodynamic equation (\ref{NonLinDiff}) given in proposition \ref{ExistenceofSolutionsTheorem} of the next section, the functions in (\ref{SufficesforGeneralInitProf}) belong in $L^\infty_{\rm{loc}}(\RR_+)$ whenever the initial profile is of class $C^{2+\theta}(\T^d;(0,\rho_c))$.

\section{Proofs}
\textbf{Proof of the One-Block estimate:}
The proof follows the one given in section 5.4 of \cite{Landim}, the main difference being at the final step of the proof, where the equivalence of ensembles proved in \cite{Stefan} and stated here as theorem \ref{EoE} is applied, instead of the one in the appendix of \cite{Landim}. As shown in \cite{Landim}, denoting 
$$V^\ell:=\bigg|\fr{(2\ell+1)^d}\sum_{|y|\leq\ell}g\big(\eta(y)\big)-\Phi\big(\eta^\ell(0)\mn\rho_c)\bigg|,$$ in order to prove the one-block estimate it suffices to prove that for all finite constants $C_0>0$ we have  
\begin{eqnarray}\label{ReduceToStatic}\qquad\lim_{\ell\ra\infty}\limsup_{N\ra\infty}\sup_{\substack{H_N(f)\leq C_0N^d\\ D_N(f)\leq C_0N^{d-2}}}\int\fr{N^d}\sum_{x\in\T_N^d}(\tau_xV^\ell)fd\nu_{\rho_*}^N\leq 0,\end{eqnarray}
where the supremum is taken among densities $f\in L^1_{+,1}(\nu_{\rho_*}^N)$.
Here, for each density $f\in L^1_{+,1}(\nu_{\rho_*}^N)$ we abbreviate by $H_N(f):=\HHH(fd\nu_{\rho_*}^N|\nu_{\rho_*}^N)$ its relative entropy with respect to $\nu_{\rho_*}^N$, and $D_N:=\mathfrak{D}_N(\sqrt{\cdot})$ where $\mathfrak{D}_N$ is the Dirichlet form associated to the generator $L_N$ of the ZRP.
Whenever $g$ is bounded, it also has sub-linear growth and a careful inspection of the proof in section 4 of \cite{Landim} shows that one can follow the steps 1 to 5 there to conclude that in order to complete the proof of the one block estimate it suffices to show that for all positive constants $C_1>0$,
\begin{eqnarray}\label{step5end}\qquad
\limsup_{\ell\ra\infty}\sup_{K\leq(2\ell+1)^dC_1}\int V^\ell(\xi)d\nu_{\ell_\star,K}(\xi)=0,\quad\ell_\star:=2\ell+1.
\end{eqnarray}
\indent The final step of the proof of the one block estimate consists in applying the equivalence of ensembles to prove (\ref{step5end}). Since the measure $\nu_{\ell_\star,K}$ is concentrated on configurations with $K$ particles, by fixing a positive integer $k$ that shall increase to infinity after $\ell$ and decomposing the set $\Lambda_\ell^d:=\big\{y\in\ZZ^d\,\big|\,|y|\leq\ell\big\}$ in cubes of side-length $k_\star=2k+1$ it follows as in \cite{Landim} that the integral $\int V^\ell d\nu_{\ell_\star,K}$ can be bounded above by
$$\int V^\ell d\nu_{\ell_\star,K}\leq
\int\bigg|\fr{k_\star^d}\sum_{|x|\leq k}g\big(\xi(x)\big)-
\Phi\Big(\frac{K}{\ell_\star^d}\Big)\bigg|d\nu_{\ell_\star,K}+C_d\frac{k\ell^{d-1}}{\ell_\star^d},$$ 
where $C_d$ is a universal constant that depends only on the dimension $d$ and on the local jump rate function $g$. Therefore, if we set
$$S(m,k):=\sup_{\substack{\ell\geq m\\K\leq\ell_\star^dC_1}}
\int\Big|\fr{k_\star^d}\sum_{|x|\leq k}g\big(\xi(x)\big)-
\Phi\Big(\frac{K}{\ell_\star^d}\Big)\Big|d\nu_{\ell_\star,K},$$
it suffices to prove that $S(m,k)$ tends to zero as $m$ and then $k$ tend to infinity. For each fixed $(m,k)\in\NN\x\NN$ we can pick a sequence $\{(\ell_n^{m,k},K^{m,k}_n)\}_{n\in\NN}$ such that $\ell_n^{m,k}\geq m$ and $K_n^{m,k}\leq(\ell_n^{m,k})_\star^dC_1$ for all $n\in\NN$, that achieves the supremum, i.e$.$ such that 
$$S(m,k)=\lim_{n\ra\infty}\int\Big|\fr{k_\star^d}\sum_{|x|\leq k}g\big(\xi(x)\big)-
\Phi\Big(\frac{K^{m,k}_n}{(\ell_n^{m,k})_\star^d}\Big)\Big|d\nu_{(\ell_n^{m,k})_\star,K^{m,k}_n}.$$
Since the sequence $\{\rho_n^{m,k}\}_{n\in\NN}$ defined by 
$\rho_n^{m,k}:=K^{m,k}_n/(\ell_n^{m,k})_\star^d$ is contained in the compact interval $[0,C_1]$, for each fixed $(m,k)\in\NN\x\NN$ we can pick a sequence $\{n_j\}_{j\in\NN}=\{n_j^{m,k}\}$ such that $\rho^{m,k}_{n_j}$ converges to some $\rho^{m,k}\in[0,C_1]$ as $j\ra\infty$.
Then since $g$ is assumed bounded it follows by the equivalence of ensembles that 
$$S(m,k)=\int\Big|\fr{k_\star^d}\sum\nolimits_{|x|\leq k}g\big(\xi(x)\big)-
\Phi\big(\rho^{m,k}\big)\Big|d\nu_{\rho^{m,k}\mn\rho_c}^{k_\star}.$$ Furthermore, for each fixed $k\in\NN$, the sequence $\{\rho^{m,k}\mn\rho_c\}_{m\in\NN}$ is contained in $[0,\rho_c]$ and thus we can choose a sequence $\{m_j\}_{j\in\NN}=\{m_j^{(k)}\}$ such that $\rho^{m_j,k}$ converges to some $\rho^k\in[0,\rho_c]$. Then by the weak continuity of the grand canonical ensemble,
\begin{eqnarray*}
\qquad\lim_{m\ra\infty}S(m,k)=\int\Big|\fr{k_\star^d}\sum\nolimits_{|x|\leq k}g\big(\xi(x)\big)-
\Phi\big(\rho^k\big)\Big|d\nu_{\rho^k}^{k_\star}
\end{eqnarray*}
for each fixed $k\in\NN$. Consequently, in order to complete the proof of the one block estimate it suffices to prove that 
\begin{eqnarray}\label{step6toend}\qquad\lim_{k\ra\infty}\int\Big|\fr{k_\star^d}\sum\nolimits_{|x|\leq k}g\big(\xi(x)\big)-
\Phi\big(\rho^k\big)\Big|d\nu_{\rho^k}^\infty=0.
\end{eqnarray}
Now, since the random variables $g\big(\eta(x)\big)$, $x\in\ZZ^d$, are uniformly bounded by $\|g\|_u$ and i.i.d$.$ with respect to $\nu_\rho^\infty$ for all $\rho\in[0,\rho_c]$, the weak law of large numbers holds in $L^2(\nu_\rho^\infty)$ uniformly over all parameters $\rho\in[0,\rho_c]$, i.e$.$ 
$$\lim_{N\ra\infty}\sup_{\rho\in[0,\rho_c]}\int\bigg|\fr{N^d}\sum_{x\in\T_N^d}g\big(\eta(x)\big)-
\Phi\big(\rho\big)\bigg|^2d\nu_\rho^\infty=0,$$ which since $\{\rho^k\}_{k\in\NN}\subs[0,\rho_c]$ proves (\ref{step6toend}) and completes the proof.$\hfill\Box$\\


\noindent\textbf{Proof of theorem \ref{MainTheorem}:} The proof follows the one in  chapter 6 of \cite{Landim}. Let $\{\mu_0^N\in\PP\MM_N^d\}$ be an initial entropy local equilibrium of profile $\rho_0$ and let $\mu_t^N:=\mu_0^NS_{tN^2}^N\in\PP\MM_N^d$, $t\geq 0$, denote the evolution of the initial distribution $\mu_0^N$ under the diffusively rescaled transition semi-group of the ZRP. Since the initial profile is continuous and takes values in $(0,\rho_c)$ there exists $\ee>0$ such that $\rho_0(\T^d)\subs(\ee,\rho_c-\ee)$. We fix such an $\ee\in(0,\rho_c/2)$. We consider the case of finite critical density. For the case $\rho_c=+\infty$ (for instance when $b\in(1,2]$ in the Evans model) one simply needs to replace $\rho_c-\ee$ by $\fr{\ee}$.\\
\indent We will first make sure that in the case of continuous sub-critical initial data the initial value problem (\ref{NonLinDiff}) with non-linearity $\Phi\in\mathcal{C}_{\rho_c}$, $\rho_c\in(0,\infty]$ admits classical solutions. This is done by using the sub-criticality of the initial data, the maximum principle and the following lemma to avoid the degeneracy of $\Phi$ at $\rho_c$. The proof of this lemma is straightforward and follows by induction.
\begin{lemma}\label{rirourin}
Let $\Phi:(0,b]\lra(0,\infty)$, $b\in(0,\infty)$, be a strictly positive $C^k$ function. There exists then large enough $M\equiv M(k)\geq 0$ such that the function $$\wt{\Phi}(\rho)=\begin{cases}\Phi(\rho),\quad& 0<\rho\leq b\\
\frac{M}{(k+1)!}(\rho-b)^{k+1}+\sum_{m=0}^k\frac{\Phi^{(m)}(b)}{m!}(\rho-b)^m,\quad&\rho\geq b\end{cases}$$ is a strictly positive $C^k$ extension $\wt{\Phi}:(0,\infty)\lra(0,\infty)$ of $\Phi$. 
\end{lemma}

 
\begin{prop}\label{ExistenceofSolutionsTheorem}
Let $\Phi\in\mathcal{C}_{\rho_c}$, $\rho_c\in(0,\infty]$ and let $\rho_0:\T^d\lra[0,\rho_c)$ be a continuous initial profile.
There exists then a unique classical solution 
$\rho\in C(\RR_+\x\T^d)$ of the initial value problem (\ref{NonLinDiff}). Furthermore, $\rho$ is $C^\infty$ on $(0,\infty)\x\T^d$ and if the initial profile is of class $C^{2+\theta}$ for some $\theta\in(0,1]$, then $\rho\in C^{1+\theta,2+\theta}(\RR_+\x\T^d)$. Moreover,
\begin{eqnarray}\label{TrulyBoundedDataAbove}\quad
\max_{u\in\T^d}\rho_t(u)\leq\max_{u\in\T^d}\rho_0(u)<\rho_c,
\end{eqnarray}
for all $t>0$. Finally, if $\rho_0$ is not constant, then for all $t>0$ we have that
\begin{eqnarray}\label{TrulyBoundedDataBelow}\quad\min_{u\in\T^d}\rho_t(u)>\min_{u\in\T^d}\rho_0(u)\geq 0.\end{eqnarray}
\end{prop}\textbf{Proof} Since $\rho_0$ is continuous and takes values in the interval $[0,\rho_c)$ it follows by the compactness of $\T^d$ that there exists $\ee>0$ such that $\max_{u\in\T^d}\rho_0(u)<\rho_c-\ee$. Then, since $\Phi '(0)\mn\Phi '(\rho_c-\ee)>0$, for any fixed $k\geq 1$ there exists by lemma \ref{rirourin} a strictly positive (two-sided) $C^k$ extension $\Psi\equiv\Psi_k:\RR\lra\RR$ of $\Phi '|_{[0,\rho_c-\ee]}$. Since $M>0$ we have that $\lim_{|\rho|\ra\infty}\Psi(\rho)=+\infty$ and therefore $c:=\inf_{\rho\in\RR}\Psi(\rho)>0$. We set $B:=\max_{\rho\in[0,\rho_c-\ee]}\Psi(\rho)$, we choose a smooth function $\chi:\RR_+\lra[0,B+1]$ such that $\chi(y)=y$ for $0\leq y\leq B$ and $\chi(y)=B+1$ for $y\geq B+1$ and consider the function $\wt{\Psi}:=\chi\circ\Psi:\RR\lra(0,\infty)$. Then its anti-derivative $\wt{\Phi}(\rho)=\int_0^\rho\wt{\Psi}(r)dr$, $\rho\in\RR$, is a $C^{k+1}$ extension of the restriction $\Phi|_{[0,\rho_c-\ee]}$ satisfying $c\leq\wt{\Phi}'(\rho)\leq B+1$ for all $\rho\in\RR$. The claim then follows by applying the results on uniformly parabolic equations obtained in \cite{Lady} (see also section 3.1.1 in \cite{Vazquez}) to the initial value problem $\pd_t\rho=\D_\Sigma\wt{\Phi}(\rho)$ with initial condition $\rho(0,\cdot)\equiv\rho_0$, for each $k\geq 1$.$\hfill\Box$\\

For the rest of the proof we denote by $\rho:\RR_+\x\T^d\lra(0,\rho_c)$ the classical solution of the initial value problem (\ref{NonLinDiff}) with $\rho(0,\cdot)\equiv\rho_0\in C^{2+\theta}(\T^d;(0,\rho_c))$, we fix $a\in(0,\rho_c)$ and denote by $\psi^{N}_t$ the Radon-Nikodym derivative of $\nu_{\rho_t(\cdot)}^N$ with respect to $\nu_a^N$, that is 
$$\psi_t^N:=\frac{d\nu_{\rho_t(\cdot)}^N}{d\nu_a^N}.$$ 
Setting $H_N(t):=H(\mu_t^N|\nu_{\rho_t()}^N)$ the relative entropy of $\mu_t^N$ with respect to $\nu_{\rho_t()}^N$ we have the following upper bound on the entropy production, proved in lemma 6.1.4 in \cite{Landim}: 
\begin{eqnarray}\label{AFirstUpBoundOnEntrProd}\qquad
\pd_tH_N(t)\leq\int\fr{\psi_t^N}\big\{N^2L_N^*\psi_t^N-\pd_t\psi_t^N\big\}d\mu_t^N,
\end{eqnarray}
for every $t\geq 0$, where $L_N^*$ is the adjoint of $L_N$ in $L^2(\nu_a^N)$. Denoting by 
\begin{eqnarray}\label{LRE}\qquad
H(t):=\limsup_{N\ra\infty}\fr{N^d}H_N(t),\quad t\in\RR_+,
\end{eqnarray} the limiting renormalised entropy, the main step in the application of the relative entropy method is to use this upper bound on $\pd_tH_N(t)$ to get an inequality of the form 
\begin{eqnarray}\label{GronwIneqOnLimitingRenormedEntropy}\qquad
H(t)\leq H(0)+ \int_0^tH(s)\beta(s)ds
\end{eqnarray} for a non-negative function $\beta$. Since $H(0)=0$ by assumption, this implies by Gronwall's inequality that $H(t)=0$ for all $t\in\RR_+$ as required. Of course in order for Gronwall's inequality to be applicable, the function $s\mapsto H(s)\beta(s)$ must belong at least in $L^1_{{\rm{loc}}}(\RR_+)$.\newpage
\begin{lemma}\label{UnifBoundNormalizedMicroscopicEntropies}
Let $\rho:\RR_+\x\T^d\lra[0,\rho_c)$ be a continuous function such that the function in (\ref{SufficesforGeneralInitProf}a) is in $L^2_{\rm{loc}}(\RR_+)$. If a sequence of initial distributions $\{\mu_0^N\}$ has relative entropy of order $o(N^d)$ with respect to $\nu_{\rho_0(\cdot)}^N$, then the upper entropy $\bbar{H}$ belongs in $L^2_{\rm{loc}}(\RR_+)$, where  
$$\bbar{H}(t):=\sup_{N\in\NN}\fr{N^d}H(\mu_t^N|\nu_{\rho_t(\cdot)}^N),\quad\;t\in\RR_+.$$ 
\end{lemma}\textbf{Proof} By remark 6.1.2 of \cite{Landim}, the relative entropy inequality shows that $\{\mu_0^N\}$ satisfies the $O(N^d)$-entropy assumption: $C(a)=\sup_{N\in\NN}\fr{N^d}H(\mu_t^N|\nu_a^N)<\infty$ for all $a\in(0,\rho_c)$. Using the relative entropy inequality once again we prove that $\bbar{H}\in L^2_{\rm{loc}}(\RR_+)$. Indeed, given $T>0$ we pick $\ee>0$ such that $\rho_c-\ee$ is an upper bound of the set $\rho([0,T]\x\T^d)$ and fix $a\in(\rho_c-\ee,\rho_c)$. By the relative entropy inequality and proposition A.1.9.1 of \cite{Landim}, according to which the function $t\mapsto H(\mu_t^N|\nu_a^N)$ is non-increasing,
\begin{eqnarray}\label{REBound2}\qquad
H_N(t)\leq\Big(1+\fr{\g}\Big)H(\mu_0^N|\nu_a^N)+\fr{\g}\log\int\Big(\frac{d\nu_a^N}{d\nu_{\rho_t(\cdot)}^N}\Big)^\g d\nu_a^N
\end{eqnarray} 
for all $t\geq 0$ and all $\g>0$. By a standard computation, 
\begin{eqnarray*}
\log\int\Big(\frac{d\nu_a^N}{d\nu_{\rho_t(\cdot)}^N}\Big)^\g d\nu_a^N&=&\sum_{x\in\T_N^d}\Big\{\g\log\frac{Z\big(\Phi(\rho_t(x/N)\big))}{Z(\Phi(a))}+
\Lambda_a\Big(\g\log\frac{\Phi(a)}{\Phi(\rho_t(x/N))}\Big)\Big\},
\end{eqnarray*}
where $\Lambda_{\rho}$ is the logarithmic moment generating function one-site ZR distribution $\nu^1_\rho$: 
$$\Lambda_\rho(r):=\log\int e^{rk}d\nu_\rho^1(k)=\log\frac{Z(e^r\Phi(\rho))}{Z(\Phi(\rho))}.$$ So if for each $t>0$ we set $$\g(t):=\fr{2}\frac{\log\frac{\phih_c}{\Phi(a)}}{\log\frac{\Phi(a)}{\Phi(m_t)}},$$ then $\g(t)\log\fr{\Phi_a(\rho_t(x/N))}\leq\fr{2}\log\frac{\phih_c}{\Phi(a)}$ for all $(t,x)\in\RR_+\x\T_N^d$, and by (\ref{REBound2}) for all $t\in[0,T]$ 
$$\bbar{H}(t)\leq\big(1+1/{\g(t)}\big)C(a)+\log Z_a(\rho_c-\ee)+\big(1/{\g(t)}\big)\log Z\big(\sqrt{\phih_c\Phi(a)}\,\big).$$
Since the function in (\ref{SufficesforGeneralInitProf}a) is in $L^2_{\rm{loc}}(\RR_+)$, the right hand side above is in $L^2([0,T])$.$\hfill\Box$\\

The bound (\ref{AFirstUpBoundOnEntrProd}) on the entropy production can be estimated explicitly. As in \cite{Landim} on one hand
\begin{eqnarray}\label{ForTaylor}\qquad
\frac{L_N^*\psi_t^N}{\psi_t^N}=
\sum_{x,y\in\T_N^d}\Big[\frac{\Phi\big(\rho_t(y/N)\big)}{\Phi\big(\rho_t(x/N)\big)}-1\Big]
\big[g\big(\eta(x)\big)-\Phi\big(\rho_t(x/N)\big)\big]p(y-x),
\end{eqnarray}
while on the other hand using the fact that $\pd_t\rho=\D_\Sigma\Phi(\rho)$ and equality (\ref{MeanDensityMainFormula}), we get
\begin{eqnarray}\label{ToAddSpatialMean}\qquad
\frac{\pd_t\psi_t^N}{\psi_t^N}=\pd_t(\log\psi_t^N)=\sum_{x\in\T_N^d}\frac{\D_\Sigma[\Phi(\rho_t)]}{\Phi(\rho_t)}\big(\frac{x}{N}\big)\Phi '\Big(\rho_t\big(\frac{x}{N}\big)\Big)\Big[\eta(x)-\rho_t\big(\frac{x}{N}\big)\Big],
\end{eqnarray}
where an integration by parts shows that in (\ref{ToAddSpatialMean}) we can replace $\eta(x)$ by $\eta^\ell(x)$, $\ell\in\ZZ_+$.\newpage
\indent Since $\Phi(\rho_t)$ is smooth and since the elementary step distribution $p$ has mean zero, we have by the second order Taylor expansion for $C^{2+\theta}$ functions that 
\begin{eqnarray*}
\sum_{y\in\T_N^d}\Big[\Phi\Big(\rho_t\big(\frac{y}{N}\big)\Big)-\Phi\Big(\rho_t\big(\frac{x}{N}\big)\Big)\Big]p(y-x)
&=&\fr{N^2}\D_\Sigma[\Phi(\rho_t)]\big(\frac{x}{N}\big)+R_t\big(\frac{x}{N}\big)\end{eqnarray*}
for all $x\in\T_N^d$, where the remainder $R_t$ satisfies
$\|R_t\|_\infty\leq
CN^{-(2+\theta)}\big\|D^2[\Phi(\rho_t)]\big\|_{\theta}$ 
for some constant $C>0$ depending only on $p\in\PP\ZZ^d$ and the dimension $d$. Therefore we can write (\ref{ForTaylor}) as
$$N^2\frac{L_N^*\psi_t^N}{\psi_t^N}
=\sum\nolimits_{x\in\T_N^d}
\Big(\frac{\D_\Sigma[\Phi(\rho_t)]}{\Phi(\rho_t)}\Big)\big(\frac{x}{N}\big)\big[g(\eta(x))-\Phi(\rho_t(x/N))\big]+r_N(t)$$ where for the remainder term $r_N(t)$ we have
$$|r_N(t)|\leq
\frac{C}{N^{\theta}}\sum_{x\in\T_N^d}\frac{\|D^2[\Phi(\rho_t)]\|_\theta}{\Phi(m_t)}\Big|g\big(\eta(x)\big)-\Phi\big(\rho_t(x/N)\big)\Big|
\leq\bar{C}N^{d-\theta}\frac{\big\|D^2[\Phi(\rho_t)]\big\|_\theta}{\Phi(m_t)},$$
with $\bar{C}:=C(\|g\|_\infty+\Phi(\rho_c-\ee))$. By this bound on the remainder and the $L^2_{\rm{loc}}(\RR_+)$-integrability of the function defined in (\ref{SufficesforGeneralInitProf}b) it follows that $\int_0^t\int r_N(s)d\mu_s^Nds=o(N^d)$ for each $t>0$.

\begin{prop}\label{AfterONEBLOCKApply} For all $t>0$ 
$$H_N(t)\leq H_N(0)+
\int_0^t\int\sum_{x\in\T_N^d}
\frac{\D_\Sigma[\Phi(\rho_s)]}{\Phi(\rho_s)}
\big(\frac{x}{N}\big)M\Big(\eta^\ell(x),\rho_s\big(\frac{x}{N}\big)\Big)d\mu_s^Nds+o_\ell(N^d),$$
where $M:\RR_+\x\RR_+\lra\RR$ is the function given by the formula 
$$M(\lambda,\rho)=\Phi(\lambda)-\Phi(\rho)-\Phi '(\rho)(\lambda-\rho)$$ and the term $o_\ell(N^d)$ satisfies $o_\ell(N^d)/N^d\lra 0$ as $N$ and then $\ell$ tend to infinity.
\end{prop}\textbf{Proof} Using the version of the one-block estimate proved in this article, the proof follows as in p. 122 of \cite{Landim} by (\ref{AFirstUpBoundOnEntrProd}) and the calculations above on the terms (\ref{ForTaylor}) and (\ref{ToAddSpatialMean}).$\hfill\Box$\\

To simplify the notation, we set $G_t:\T^d\x\RR_+\lra\RR$, $t\geq 0$, the function defined by 
\begin{eqnarray}\label{G}\qquad
G_t(u,\lambda)=\frac{\D_\Sigma[\Phi(\rho_t)]}{\Phi(\rho_t)}(u)M\big(\lambda,\rho_t(u)\big).
\end{eqnarray}
By the relative entropy inequality we get that 
\begin{eqnarray*}
\int\sum_{x\in\T_N^d}
G_s\Big(\frac{x}{N},\eta^\ell(x)\Big)d\mu_s^N&\leq&\fr{\g_s}H_N(s)+\fr{\g_s}\log\int e^{\g_s\sum_{x\in\T_N^d}G_s(\frac{x}{N},\eta^\ell(x))}d\nu_{\rho_s(\cdot)}^N
\end{eqnarray*}
for any positive function $(0,\infty)\ni s\mapsto\g_s\in(0,\infty)$ and each $s>0$. We combine this inequality with proposition \ref{AfterONEBLOCKApply}, divide by $N^d$ and take $\limsup$ first as $N\ra\infty$ and then as $\ell\ra\infty$. Then if the function $\g$ can be chosen so that $\beta:=1/\g\in L^2_{\rm{loc}}(\RR_+)$, we can use lemma \ref{UnifBoundNormalizedMicroscopicEntropies} to pass the $\limsup$ as $N\ra\infty$ inside the time integral of $s\mapsto H_N(s)\beta(s)/N^d$ to get Gronwall's inequality (\ref{GronwIneqOnLimitingRenormedEntropy}) but with the term 
\begin{eqnarray}\label{GronPlusTerm}\qquad \limsup_{\ell\ra\infty}\limsup_{N\ra\infty}\fr{ N^d}\int_0^t\fr{\g_s} \log\int e^{\g_s\sum_{x\in\T_N^d}G_s(\frac{x}{N},\eta^\ell(x))}d\nu_{\rho_s(\cdot)}^Nds\end{eqnarray}
added to its right hand side.\\
\indent So the rest of the proof is devoted to proving that the function $\beta\equiv 1/\g\in L^2_{\rm{loc}}(\RR_+)$ can be chosen so that for each time $t>0$ the term in (\ref{GronPlusTerm}) is non-positive. We begin by noting that the function $G:\RR_+\x\T^d\x\RR_+\lra\RR$ defined in (\ref{G}) satisfies the inequality
\begin{eqnarray}\label{UnifBoundOnGOnFinTImeHorizon}\qquad
\sup\nolimits_{u\in\T^d}|G_t(u,\lambda)|\leq C\cdot C_t\cdot(1+\lambda) \quad\mbox{for all }t,\lambda>0
\end{eqnarray} where
$C=\big\{\Phi(\rho_c-\ee)+\max_{r\in[0,\rho_c-\ee]} r\Phi '(r)\big\}\mx 2\|g'\|_\infty<+\infty$, $C_t:=\big\|{\Delta_{\Sigma}\Phi(\rho_t)}/{\Phi(\rho_t)}\big\|_\infty$. For each $K>1$ we denote by $\g^K:(0,\infty)\lra(0,\infty)$ the function $\g_t^K:=\fr{KCC_t}\log\frac{\phih_c}{\Phi(\rho_c-\ee)}$. Since the function in (\ref{SufficesforGeneralInitProf}b) is in $L^2_{\rm{loc}}(\RR_+)$, the function $\beta^K:=1/\g^K$ belongs in $L^2_{\rm{loc}}(\RR_+)$. Using inequality (\ref{UnifBoundOnGOnFinTImeHorizon}) and the $L_{\rm{loc}}^2(\RR_+)$-integrability of $\beta^K$ it is straightforward to check that the family $\{h_K^{N,\ell}\}_{(N,\ell)\in\NN}$ of the functions 
$$h_K^{N,\ell}(t)=\fr{\g^K_tN^d}\log\int e^{\g^K_t\sum_{x\in\T_N^d}G_t(\frac{x}{N},\eta^\ell(x))}d\nu_{\rho_t(\cdot)}^N,\quad t\geq 0$$ 
is dominated by an $L_{\rm{loc}}^2(\RR_+)$-function for each $K>1$. This permits to pass the superior limits as $N\ra\infty$ and then $\ell\ra\infty$ inside the time integral in (\ref{GronPlusTerm}) for each $K>1$. Consequently, in order to complete the proof it suffices to show that we can choose $K>1$ so that for each $t>0$, \begin{eqnarray}\label{ToConcludeTheProofOfTheRelEntrIneq2}\qquad
\limsup_{\ell\ra\infty}\limsup_{N\ra\infty}\fr{N^d}\log\int e^{\g^K_t\sum_{x\in\T_N^d}G_t(\frac{x}{N},\eta^\ell(x))}d\nu_{\rho_t(\cdot)}^N\leq 0.\end{eqnarray} 
\indent This will follow from the estimate of the following lemma. Before we state this lemma, we mention that by Cramer's theorem, for each $\rho\in[0,\rho_c)$ the family $\{\eta(x)\}_{x\in\T_N^d}$ satisfies with respect to $\nu_\rho^\infty$ the Large Deviations Principle (LDP) with speed $N^d$ and rate functional the Legendre transform 
\begin{eqnarray*}
\Lambda_{\rho}^*(\lambda):=\sup_{r\in\RR}\big\{\lambda r-\Lambda_\rho(r)\big\}=
\begin{cases}\lambda\log\frac{\Phi(\lambda\mn\rho_c)}{\Phi(\rho)}-\log\frac{Z(\Phi(\lambda\mn\rho_c))}{Z(\Phi(\rho))}\quad&{\rm{if }}\;\lambda\geq 0\\
+\infty\quad&{\rm{if }}\;\lambda<0
\end{cases}
\end{eqnarray*}

\begin{lemma}\label{LDPlemma}  Let $\rho:\T^d\lra(0,\rho_c-\ee)$, $\ee\in(0,\rho_c)$, be a continuous profile and let $G:\T^d\x\RR_+\lra\RR$ be a continuous function such that $\sup_{u\in\T^d}G(u,\lambda)\leq C_0+C_1\lambda$ holds for some constants $C_0\geq 0$ and $C_1\in[0,\fr{2}\log\frac{\phih_c}{\Phi(\rho_c-\ee)})$. Then 
$$\limsup_{\ell\ra\infty}\limsup_{N\ra\infty}\fr{N^d}\log\int e^{\sum_{x\in\T_N^d}G(\frac{x}{N},\eta^\ell(x))}d\nu_{\rho(\cdot)}^N\leq
\int_{\T^d}\sup_{\lambda>0}\Big\{G(u,\lambda)-\fr{2}\Lambda_{\rho(u)}^*(\lambda)\Big\}du.$$
\end{lemma}\textbf{Proof} The proof is based on the LDP described above and the Laplace-Varadhan lemma. It is given in lemma 6.1.8 of \cite{Landim} in the case $\rho_c=\infty$. A careful inspection of the proof there shows that it also applies to the case $\rho_c<\infty$.$\hfill\Box$\\

 We recall that $G:[0,T]\x\T^d\x\RR_+\lra\RR$ satisfies the bound (\ref{UnifBoundOnGOnFinTImeHorizon}). Therefore if we choose $K>2$ then the function $\g_t^K G_t$ satisfies the assumptions of lemma \ref{LDPlemma} for each fixed $t>0$, and so for each $K>2$ the term in (\ref{ToConcludeTheProofOfTheRelEntrIneq2}) is bounded above by $$\int_{\T^d}\sup_{\lambda>0}\Big\{\g_t^K G_t(u,\lambda)-\fr{2}\Lambda^*_{\rho_t(u)}(\lambda)\Big\}du.$$
To complete the application of the relative entropy method it remains to show that, by enlarging $K>2$ if necessary, this last term is non-positive for all $t>0$. This follows by the next lemma, which is a simple generalisation of lemma 6.1.10 in \cite{Landim}.
\begin{lemma}
For every $\ee>0$,
$$\sup_{(\lambda,\rho)\in(0,\infty)\x(0,\rho_c-\ee]}\frac{|M(\lambda,\rho)|}{\Lambda_\rho^*(\lambda)}<+\infty.$$
\end{lemma}\textbf{Proof} We first choose $\delta\in(0,\frac{\ee}{2})$. We decompose the set $(0,\infty)\x(0,\rho_c-\ee]$ in two disjoint subsets ($\lambda\backsim\rho$ and $\lambda\gg\rho$) and prove the claim on each. We start with the region $\lambda\backsim\rho$:
$$\mathcal{E}_1:=
\big\{(\lambda,\rho)\in\RR_+\x(0,\rho_c-\ee]\,\big|\,0<\lambda\leq\rho_c-\ee+\delta\big\},$$ where we recall that if $\rho_c=+\infty$, $\rho_c-\ee$ is to be interpreted as $1/\ee$. By the Taylor expansion of $\Phi$ around the point $\rho\in(0,\rho_c)$, we have
$M(\lambda,\rho)=\int_\rho^\lambda\Phi ''(r)(\lambda-r)dr$ for all $\lambda,\rho\in(0,\rho_c)$. So since $\{\lambda|(\lambda,\rho)\in\mathcal{E}_1\;\mbox{for some }\rho\in(0,\rho_c-\ee]\}\subs(0,\rho_c)$, $$|M(\lambda,\rho)|\leq\frac{A_1}{2}(\lambda-\rho)^2\quad\mbox{for all }(\lambda,\rho)\in\mathcal{E}_1$$
where $A_1:=\sup_{0\leq r\leq\rho_c-\ee+\delta}|\Phi ''(r)|<+\infty$. For the denominator we note that the rate functional $\Lambda_\rho^*$ is $C^1$ on $(0,\infty)$ and $C^2$ on $(0,\rho_c)$ with 
$$\frac{d}{d\lambda}\Lambda_\rho^*(\lambda)=\log\frac{\Phi(\lambda\mn\rho_c)}{\Phi(\rho)},\quad\lambda>0,
\qquad\frac{d^2}{d\lambda^2}\Lambda_\rho^*(\lambda)=\frac{\Phi '(\lambda)}{\Phi(\lambda)},\quad\lambda\in(0,\rho_c).$$
Since $\Lambda_\rho^*$ and its derivative vanish at $\rho$, by the Taylor expansion of $\Lambda_\rho^*$ around $\rho\in(0,\rho_c)$ we have that 
$\Lambda_\rho^*(\lambda)=\int_\rho^\lambda(\Lambda_\rho^*)''(r)(\lambda-r)dr$
for all $\lambda\in(0,\rho_c)$ and therefore  
$$\Lambda_\rho^*(\lambda)\geq \frac{B_1}{2}(\lambda-\rho)^2\quad\mbox{for all }(\lambda,\rho)\in\mathcal{E}_1,$$ where 
$B_1:=\inf_{0<r\leq\rho_c-\ee+\delta}(\Lambda_\rho^*)''(r)>0$. Combining these estimates, we get the required bound on the region $\mathcal{E}_1$.\\ 
\indent We turn now to the set 
$$\mathcal{E}_2=\big\{(\lambda,\rho)\in(0,\infty)\x(0,\rho_c-\ee]\,\big|\,\lambda>\rho_c-\ee+\delta\big\}.$$
Note that for all $(\lambda,\rho)\in\mathcal{E}_2$ we have that $\lambda>\rho+\delta$. Recalling that $\Phi$ is Lipschitz with Lipschitz constant $\leq\|g'\|_u$, we get an upper bound for the numerator
$$|M(\lambda,\rho)|\leq  2\|g'\|_\infty\lambda \quad \mbox{for all }(\lambda,\rho)\in\mathcal{E}_2.$$
Since $\Lambda_\rho^*$ is convex as the supremum of linear functions we also have
$$\Lambda_\rho^*(\lambda)\geq A_2+B_2\cdot(\lambda-\rho_c-\ee+\delta) \quad \mbox{for all }(\lambda,\rho)\in\mathcal{E}_2, $$
where $A_2=\inf_{\rho\in(0,\rho_c-\ee]}\Lambda_\rho^*(\rho_c-\ee+\delta)>0$ and $B_2=\inf_{\rho\in(0,\rho_c-\ee]}(\Lambda_{\rho}^*)'(\rho_c-\ee+\delta)>0$. The last two displays together imply the required bound on the region $\mathcal{E}_2$. This completes the proof of the lemma and the application of the relative entropy method.$\hfill\Box$\\

\hspace{-14.5pt}
\noindent{\footnotesize{{\bf{Acknowledgements:}} The current research has been co-financed by the European Union (European Social Fund- ESD) and by national resources through the operational programme ``Education and Lifelong
Learning" of the National Strategic Research Frame (NSRF) - Financed Research Project: Herakleitos II. Investment in the society of knowledge through the European Social Fund.}}\\
\noindent{\footnotesize{\indent I am thankful to my advisor Michail Loulakis for his continuous help and encouragement and to Claudio Landim for useful conversations.}} 
\bibliographystyle{plain}
\bibliography{paperbib}
\end{document}